\documentstyle[12pt,amssymb,amstex]{amsart}

\def\cs{{\cal S}}

\def\ga{{\frak A}}
\def\gb{{\frak B}}

\def\bc{{\mathbb C}}

\def\bn{{\mathbb N}}

\def\a{\alpha}

\def\l{\lambda}

\def\p{\psi}

\def\f{\varphi} 
\def\v{\phi}
  
\def\w{\omega} 

\def\xb{{\mathbf{x}}}
\def\yb{{\mathbf{y}}}

\newtheorem{thm}{Theorem}[section]

\newtheorem{cor}[thm]{Corollary}
\newtheorem{prop}[thm]{Proposition}
\newtheorem{defin}[thm]{Definition}

\theoremstyle{remark}
\newtheorem{ex}{Example}[section]

\def\id{{\bf 1}\!\!{\rm I}}
\begin{document}

\title[On tensor products of weak mixing sequences]
{On tensor products of weak mixing vector sequences and their
applications to uniquely $E$-weak mixing $C^*$- dynamical systems}
\author{Farrukh Mukhamedov}
\address{Farrukh Mukhamedov\\
 Department of Computational \& Theoretical Sciences\\
Faculty of Science, International Islamic University Malaysia\\
P.O. Box, 141, 25710, Kuantan\\
Pahang, Malaysia} \email{{\tt far75m@@yandex.ru}}

\begin{abstract}
We prove that, under certain conditions, uniform weak mixing (to
zero) of the bounded sequences in Banach space implies uniform weak
mixing of its tensor product. Moreover, we prove that ergodicity of
tensor product of the sequences in Banach space implies its weak
mixing. Applications of the obtained results, we prove that tensor
product of uniquely $E$-weak mixing $C^*$-dynamical systems is also
uniquely $E$-weak mixing as well.

 \vskip 0.3cm \noindent {\it Mathematics Subject
Classification}: 46L35, 46L55, 46L51, 28D05
60J99.\\
{\it Key words}: Uniform weak mixing, weak mixing, ergodicity,
uniquely E-weak mixing, $C^*$-dynamical system.
\end{abstract}

\maketitle


\section{Introduction}

Let $X$ be a Banach spaces with dual space $X^*$. In what follows
$B_X$ denotes the unit ball in $X$, i.e. $B_X=\{x\in X: \|x\|\leq
1\}$.

Recall that a sequence $\{x_k\}$ in $X$ is said to be
\begin{enumerate}

\item[(i)] {\it weakly mixing to zero} if

$$\lim\limits_{n\to\infty}\frac1n\sum\limits_{k=1}^n|f(x_k)|=0, \
\ \ \textrm{for all} \ \ f\in X^*;$$

\item[(ii)] {\it uniformly weakly mixing to zero} if
$$
\lim\limits_{n\to\infty}\sup\bigg\{\frac1n\sum\limits_{k=1}^n|f(x_k)|
:  \ f\in B_{X^*}\bigg\}=0;
$$

\item[(iii)] {\it weakly ergodic} if
$$
\lim\limits_{n\to\infty}\frac1n\bigg|\sum\limits_{k=0}^{n-1}f(x_k)\bigg|=0
\ \ \textrm{for all} \ \ f\in X^*;$$

\item[(iii)] {\it ergodic} if
$$
\lim\limits_{n\to\infty}\frac1n\bigg\|\sum\limits_{k=0}^{n-1}x_k\bigg\|=0.
$$
\end{enumerate}

From the definitions one can see that uniform weakly mixing
implies weakly mixing, as well as ergodicity implies weak
ergodicity. But, the converse is not true at all.

\begin{ex} \cite{Z} Let  $X=L^2([0,1])$ and
$1<n_1<n_2\cdots$ be a sequence in $\bn$ such that
$$\frac{n_j-1}{n_{j+1}-1}\le \frac12, \ \ \ j\in\bn
$$
(for example, $n_1=1, \ n_2=2$ and $n_{j+1}=2n_j-1$ for
$j\in\bn$).  Let
$$1>t_1>t_2>\cdots >0, \ \ t_j\rightarrow 0$$ be real numbers and
$g_i:[0,1]\rightarrow [0,\infty), \ j\in\bn$ be continuous
functions such that
$$
supp(g_i)\subset [t_{j+1},t_j] \ \ \textrm{and} \ \ \|g_j\|_2=1
$$ for all $j\in\bn.$

Put
$$
f_k=g_i \ \ \textrm{for} \ \  n_j\le k\le n_{j+1},$$ then
$(f_k)_{k\geq 1}$ is a bounded sequence in $L^2([0,1]),$ which is
weakly convergent to zero, and so is weakly mixing to zero, but
which is not uniformly weakly mixing to zero.

\end{ex}

Recall \cite{Z} that a sequence $\{x_k\}$ in a Banach space $X$ is
called {\it convex shift-bounded} if there exists a constant $c>0$
such that
\begin{equation*}
\bigg\|\sum_{j=1}^p\l_jx_{j+k}\bigg\|\leq c
\bigg\|\sum_{j=1}^p\l_jx_{j}\bigg\|, \ \ k\geq 1
\end{equation*}
holds for any $p\in\bn$ and $\l_1,\dotsm\l_p\geq 0$.

One can see that every convex shift-bounded sequence is bounded.

\begin{ex}\label{ex2} Let $U:X\to X$ be a power bounded linear operator (i.e. the sequence
$\{\|U^k\|\}$ is bounded). Take $x\in X$ then the sequence
$\{U^k(x)\}$ is convex shift-bounded.
\end{ex}

The following theorem (see \cite{Z}) characterizes weak mixing to
zero which is a counter part of the Blum-Hanson theorem
\cite{BH},\cite{JL}.

\begin{thm}\label{Zi0} For a convex shift-bounded sequence $\{x_k\}$ in a Banach space $X$ the
following conditions are equivalent:
\begin{itemize}
\item[(i)] $\{x_k\}$ is  weakly mixing to zero; \item[(ii)]
$\{x_k\}$ is uniformly weakly mixing to zero;
\end{itemize}
\end{thm}

There is also a characterization of uniformly weak mixing to zero
by mean egodic convergence.

\begin{thm}\label{Zi}
For a bounded sequence $\{x_k\}$ in a Banach space $X$ the
following conditions are equivalent:
\begin{itemize}
\item[(i)] $\{x_k\}$ is uniformly weakly mixing (resp. weakly
mixing) to zero;

\item[(ii)] For every sequence $k_1<k_2,\cdots$ in $\bn$ with
$\sup\limits_{n\in\bn}\frac{k_n}n<+\infty$ the sequence
$\{x_{k_n}\}$ is ergodic (resp. weakly ergodic).
\end{itemize}
\end{thm}

From this theorem we conclude that weakly ergodicity does not
imply ergodicity too.

In the mentioned and others related papers (see \cite{BB, JL, K})
tensor product of sequences which obey mixing and ergodicity were
not considered. Section 2 of this note is devoted to the extension
of the well-known classical results, stating that a transformation
is weakly mixing if and only if its Cartesian square is ergodic
\cite{ALW}, for the tensor product of sequences in Banach spaces. In
next section 3, we provide some applications of the obtained results
to uniquely $E$-ergodic, uniquely $E$-weak mixing $C^*$-dynamical
systems. Note that such dynamical systems were investigated in
\cite{AD,FM,FM1,M2,MT}.

\section{Weak mixing vector sequences}

Let $X$, $Y$ be two Banach spaces with dual spaces $X^*$ and $Y^*$,
respectively. Completion of the algebraic tensor product $X\odot Y$
with respect to a cross norm $\alpha$ is denoted by $X\otimes_\alpha
Y$. By $\alpha^*$ we denote conjugate cross norm to $\alpha$ defined
on $X^*\odot Y^*$.

For the dual Banach spaces $X^*$ and $Y^*$ denote
\begin{eqnarray*}
B_{X^*}\odot B_{Y^*}=\bigg\{\sum_{k=1}^n\lambda_k x_k\otimes y_k
\bigg|&& \ \{x_k\}_{k=1}^n\subset B_{X^*},\ \{y_k\}_{k=1}^n\subset
B_{Y^*},\\
&& \ \lambda_k\geq 0,\  \sum_{k=1}^n\lambda_k\leq 1, \ n\in\bn
\bigg\}.
\end{eqnarray*}

By $B_{X^*}\otimes_{\a^*} B_{Y^*}$ denote the closure of
$B_{X^*}\odot B_{Y^*}$ with respect to conjugate cross-norm
$\a^*$. One can see that $B_{X^*}\otimes_{\a^*} B_{Y^*}\subset
B_{(X\otimes_{\alpha}Y)^*}$.

In what follows we consider the following two conditions:

\begin{enumerate}
\item[(I)] $B_{X^*}\otimes_{\a^*} B_{Y^*}=B_{(X\otimes_{\alpha}Y)^*}$.

 \item[(II)] $X^*\otimes_{\alpha^*}Y^*=(X\otimes_\alpha Y)^*$.

\end{enumerate}

One has the following

\begin{prop}\label{conj} Let $X$ and $Y$ be Banach spaces with a cross-norm $\a$
such that the property (I) holds. Then (II) is satisfied.
\end{prop}

\begin{pf} Assume that (I) is satisfied.  Now let us take an arbitrary
$f\in(X\otimes_\a Y)^*$, and show that it can be approximated by
elements of $X^*\otimes_{\a^*}Y^*$. Indeed, denote
$g=\frac{f}{\|f\|}$. Then $g\in B_{(X\otimes_\a Y)^*}$. Due to (I)
we conclude that $g\in X^*\otimes_{\a^*} Y^*$. Hence, $f=\|f\|g$
belongs to $X^*\otimes_{\a^*} Y^*$.
\end{pf}

In what follows, for given $r>0$ and $a\in X$ denote
$$
B_{r,X}(a)=\{x\in X:\ \|x-a\|\leq r\}.
$$

\begin{prop}\label{conj1} Let $X$ and $Y$ be Banach spaces with a cross-norm $\a$.
Then the property
(I) is satisfied if and only if there is a number $r>0$ ($r\leq 1$)
and an element $y\in X^*\otimes_{\a^*} Y^*$ such that
\begin{equation}\label{Bry}
B_{r,(X\otimes_\a Y)^*}(y)\subset B_{X^*}\otimes_{\a^*} B_{Y^*};
\end{equation}.
\end{prop}

\begin{pf}  It is evident that (I) implies the last property, since it is satisfied with
$r=1$ and $y=0$. Now prove the reverse implication, i.e. assume that
there is $r_0>0$ and an element $y_0\in X^*\otimes_{\a^*} Y^*$ such
that \eqref{Bry} holds. We readily see that $y_0\in
B_{X^*}\otimes_{\a^*} B_{Y^*}$. To prove the statement, it is enough
to establish that $B_{(X\otimes_{\alpha}Y)^*}\subset
B_{X^*}\otimes_{\a^*} B_{Y^*}$. Take any $x\in
B_{(X\otimes_{\alpha}Y)^*}$. Consider an element $z=y_0+r_0x$, which
clearly belongs to $B_{r_0,(X\otimes_{\a}Y)^*}$. Due to the
assumption, we conclude that $z\in B_{X^*}\otimes_{\a^*} B_{Y^*}$,
therefore, one gets that $x=\frac{1}{r_0}(z-y_0)$ belongs to
$B_{X^*}\otimes_{\a^*} B_{Y^*}$.
\end{pf}

\begin{ex} Let us give some more example which satisfy (I)
and (II) conditions.
\begin{enumerate}
\item[(i)] Let $1<p,q<\infty$, with conjugate indices $p',q'$ (i.e. $p'=\frac{p}{p-1}$).
Consider  $\ell_p$, $\ell_q$. Then for the projective norm $\pi$ one
has $(\ell_p\otimes_\pi\ell_q)^*=\ell_{p'}\otimes_{\pi^*}\ell_{q'}$
if and only if $p>q'$ (see Corollary 4.24, Theorem 4.21 \cite{R}).

\item[(ii)] We give here a sufficient condition to satisfy (II). The
proof can be found in (see Theorem 5.33 \cite{R}).

Let $X$ and $Y$ be Banach spaces such that $X^*$ has the
Radon-Nikodym property and either $X^* $ or $Y^*$ has the
approximation property. Then
$$
(X\otimes_\epsilon Y)^*=X^*\otimes_\pi Y^*
$$
here $\epsilon$ and $\pi$ are the injective and the projective
norms, respectively.

Note that more examples can be found in \cite{R}.
\end{enumerate}
\end{ex}

\begin{thm}\label{UWM}
Let $X$ and $Y$ be two Banach spaces with a cross-norm $\a$ such
that the property (I) is satisfied. Let $\{x_k\}$ be a bounded
sequence in $X$. Then the following assertions are equivalent
\begin{itemize}
\item[(i)] for any bounded sequence $\{y_k\}$ in $Y$, the sequence
$\{x_k\otimes y_k\}$ in $X\otimes_\a Y$ is uniformly weakly mixing
to zero; \item[(ii)] $\{x_k\}$ is uniformly weakly mixing to zero.
\end{itemize}
\end{thm}
\begin{pf} (i)$\Rightarrow$ (ii). Let us take any nonzero element
$y\in Y$. Define a sequence $\{y_k\}$ by $y_k=y$ for all
$k\in\bn$. For the defined sequence due to condition (i) we have
\begin{equation}\label{U1}
\lim\limits_{n\to\infty}\sup\bigg\{\frac1n\sum\limits_{k=1}^n|f(x_k\otimes
y)|: f\in B_{(X\otimes_{\alpha}Y)^*}\bigg\}=0.
\end{equation}
Now take $f=g\otimes h$ with $g\in B_{X^*}$ and $h\in B_{Y^*}$,
$h(y)\neq 0$. Then from \eqref{U1} one gets
\begin{equation*}
\lim\limits_{n\to\infty}\bigg(\sup_{g\in
B_{X^*}}\bigg\{\frac1n\sum\limits_{k=1}^n|g(x_k)|\bigg\}\bigg)|h(y)|=0
\end{equation*}
which implies the assertion.

(ii)$\Rightarrow$ (i).  Let $\{y_k\}$ be an arbitrary bounded
sequence in $Y$, and $f\in B_{X^*}$, $g\in B_{Y^*}$ be any
functionals. Then the Schwarz inequality yields
\begin{eqnarray}\label{U2}
\frac1n\sum\limits_{k=1}^n|f(x_k)g(y_k)|&\leq &
\sqrt{\frac1n\sum\limits_{k=1}^n|f(x_k)|^2}\sqrt{\frac1n\sum\limits_{k=1}^n|g(y_k)|^2}\nonumber\\
&\leq&
\max_k\{\|y_k\|\}\|g\|\sqrt{\frac1n\sum\limits_{k=1}^n|f(x_k)|^2}.
\end{eqnarray}
Moreover,
\begin{equation*}
\sup_{f\in
B_{X^*}}\bigg\{\frac1n\sum\limits_{k=1}^n|f(x_k)|^2\bigg\}\leq
\max\{\|x_k\|\}\sup_{f\in
B_{X^*}}\bigg\{\frac1n\sum\limits_{k=1}^n|f(x_k)|\bigg\}\longrightarrow
0 \ \ \textrm{as} \ \ n\to\infty.
\end{equation*}
Therefore, \eqref{U2} implies that
\begin{eqnarray}\label{U3}
\lim_{n\to\infty}\sup_{f\in B_{X^*},\atop g\in
B_{Y^*}}\bigg\{\frac1n\sum\limits_{k=1}^n|f\otimes g(x_k\otimes
y_k)|\bigg\}=0.
\end{eqnarray}

Hence, using the norm-denseness of the elements $\sum_{k=1}^m\l_k
f_k\otimes g_k$, $\{f_k\}\subset B_{X^*}$, $\{g_k\}\subset B_{Y^*}$
(where $\l_k\geq 0$, $\sum_{k=1}^n\l_k\leq 1$ ) in
$B_{X^*}\otimes_{\a^*}B_{Y^*}$ from \eqref{U3} one gets
\begin{eqnarray}\label{U4}
\lim_{n\to\infty}\sup_{\f\in
B_{X^*}\otimes_{\a^*}B_{Y^*}}\bigg\{\frac1n\sum\limits_{k=1}^n|\f(x_k\otimes
y_k)|\bigg\}=0.
\end{eqnarray}

Thanks to property (I) one has
\begin{equation*}
\sup_{f\in B_{(X\otimes_\a
Y)^*}}\bigg\{\frac{1}{n}\sum_{k=0}^{n-1}|f(x_k\otimes y_k)|\bigg\}=
\sup_{w\in
B_{X^*}\otimes_{\a^*}B_{Y^*}}\bigg\{\frac{1}{n}\sum_{k=0}^{n-1}|w(x_k\otimes
y_k)|\bigg\},
\end{equation*}
consequently \eqref{U4} yields the required statement.
\end{pf}

{\bf Remark.} From the proof of Theorem \ref{UWM} one can see that
the implication (i)$\Rightarrow$ (ii) is still valid without
property (I).

Using the same argument as above given the proof we get the
following

\begin{thm}\label{WM}
Let $X$ and $Y$ be two Banach spaces with a cross-norm $\a$ such
that property (II) is satisfied. Let $\{x_k\}$ be a bounded sequence
in $X$. Then the following assertions are equivalent
\begin{itemize}
\item[(i)] for any bounded sequence $\{y_k\}$ in $Y$, the sequence
$\{x_k\otimes y_k\}$ in $X\otimes_\a Y$ is weakly mixing to zero;
\item[(ii)] $\{x_k\}$ is weakly mixing to zero.
\end{itemize}
\end{thm}

\begin{prop}\label{EU}
Let $X$ be a Banach space and $\{x_k\}$ be a bounded sequence in
$X$ such that the sequence $\{x_k\otimes x_k\}$ is ergodic in
$X\otimes _\a X$. Then $\{x_k\}$ is uniformly weakly mixing to
zero.
\end{prop}

\begin{pf} Ergodicity of the the sequence $\{x_k\otimes x_k\}$ means that
\begin{equation}\label{xx0}
\lim\limits_{n\to\infty}\frac1n\bigg\|\sum\limits_{k=1}^{n}x_k\otimes
x_k\bigg\|=0.
\end{equation}
Due to equality
$$
\sup_{{\mathbf{f}}\in B_{(X\otimes_\a
Y)^*}}\bigg|{\mathbf{f}}\bigg(\frac{1}{n}\sum\limits_{k=0}^{n-1}x_k\otimes
x_k\bigg)\bigg|=\frac1n\bigg\|\sum\limits_{k=1}^{n}x_k\otimes
x_k\bigg\|
$$
one finds
\begin{equation}\label{xx1}
\sup_{f\in B_{X^*}}\bigg\{\frac{1}{n}\bigg|f\otimes
f\bigg(\sum\limits_{k=0}^{n-1}x_k\otimes x_k\bigg)\bigg|\bigg\}\leq
\frac1n\bigg\|\sum\limits_{k=1}^{n}x_k\otimes x_k\bigg\|.
\end{equation}

On the other hand, we have
\begin{eqnarray*}\label{EU1}
\sup_{f\in B_{X^*}}\bigg\{\frac{1}{n}\bigg|f\otimes
f\bigg(\sum\limits_{k=0}^{n-1}x_k\otimes x_k\bigg)\bigg|\bigg\}&=&
\sup_{f\in
B_{X^*}}\bigg\{\frac{1}{n}\bigg|\sum\limits_{k=0}^{n-1}f\otimes
f(x_k\otimes x_k)\bigg|\bigg\}\nonumber \\
&=& \sup_{f\in
B_{X^*}}\bigg\{\frac{1}{n}\sum\limits_{k=0}^{n-1}|f(x_k)|^2\bigg\}
\end{eqnarray*}
which with \eqref{xx0},\eqref{xx1} yields
$$
 \lim\limits_{n\to\infty}\sup_{f\in
B_{X^*}}\bigg\{\frac{1}{n}\sum\limits_{k=0}^{n-1}|f(x_k)|^2\bigg\}=0.
$$

Hence, the Schwarz inequality implies that
$$
\sup_{f\in
B_{X^*}}\bigg\{\frac{1}{n}\sum\limits_{k=0}^{n-1}|f(x_k)|\bigg\}\leq
\sqrt{\sup_{f\in
B_{X^*}}\bigg\{\frac{1}{n}\sum\limits_{k=0}^{n-1}|f(x_k)|^2\bigg\}}
$$
Therefore, we find that $\{x_k\}$ is uniformly weakly mixing to
zero.
\end{pf}

Similarly, one can prove

\begin{prop}\label{WEW}
Let $X$ be a Banach space and $\{x_k\}$ be a bounded sequence in
$X$ such that the sequence $\{x_k\otimes x_k\}$ is weakly ergodic
in $X\otimes _\a X$. Then $\{x_k\}$ is weakly mixing to zero.
\end{prop}

\begin{thm}\label{EUWM}
Let $X$ be a Banach spaces with a cross-norm $\a$ on $X\odot X$
such that condition (I) is satisfied with $Y=X$. Let $\{x_k\}$ be
a be bounded sequence in $X$. Then the following assertions are
equivalent
\begin{itemize}
\item[(i)] the sequence $\{x_k\otimes x_k\}$ is ergodic in
$X\otimes _\a X$; \item[(ii)] the sequence $\{x_k\otimes x_k\}$ is
uniformly weakly mixing to zero in $X\otimes _\a X$; \item[(iii)]
$\{x_k\}$ is uniformly weakly mixing to zero.
\end{itemize}
\end{thm}

\begin{pf} The implication (i)$\Rightarrow$ (iii) immediately
follows from Proposition \ref{EU}. The implication
(iii)$\Rightarrow$ (ii) follows from Theorem \ref{UWM}. The
implication (ii)$\Rightarrow$ (i) is evident.
\end{pf}

Using the same argument as above given the proof with Theorem
\ref{WM} one gets the following

\begin{thm}\label{EWM}
Let $X$ be a Banach spaces with a cross-norm $\a$ on $X\odot X$ such
that condition (II) is satisfied with $Y=X$. Let $\{x_k\}$ be a
bounded sequence in $X$. Then the following assertions are
equivalent
\begin{itemize}
\item[(i)] the sequence $\{x_k\otimes x_k\}$ is weakly ergodic in
$X\otimes _\a X$; \item[(ii)] the sequence $\{x_k\otimes x_k\}$ is
weakly mixing to zero in $X\otimes _\a X$; \item[(iii)] $\{x_k\}$
is weakly mixing to zero.
\end{itemize}
\end{thm}

\begin{thm}\label{UWMM}
Let $X$ and $Y$ be two Banach spaces with a cross-norm $\a$ on
$X\odot Y$ such that condition (I) (resp. (II)) is satisfied. Let
$\{x_k\}$ be a bounded sequence in $X$. The following assertions
are equivalent
\begin{itemize}
\item[(i)] for any bounded sequence $\{y_k\}$ in $Y$, the sequence
$\{x_k\otimes y_k\}$ in $X\otimes_\a Y$ is ergodic (resp. weakly
ergodic); \item[(ii)] $\{x_k\}$ is uniformly weakly mixing (resp.
weakly mixing) to zero.
\end{itemize}
\end{thm}

\begin{pf} (i)$\Rightarrow$ (ii). Let us take any nonzero element
$y\in Y$. Define a sequence $\{y_k\}$ by $y_k=y$ for all
$k\in\bn$. For the defined sequence due to condition (i) we have
\begin{equation}\label{U5}
\lim\limits_{n\to\infty}\bigg\|\frac1n\sum\limits_{k=1}^nx_k\otimes
y\bigg\|=\lim\limits_{n\to\infty}\bigg\|\frac1n\sum\limits_{k=1}^nx_k\bigg\|\|y\|=0
\end{equation}
which means $\{x_k\}$ is ergodic. The condition yields that
$\{x_k\otimes x_k\}$ is ergodic, hence Theorem \ref{EUWM} implies
that that $\{x_k\}$ is uniformly weakly mixing to zero.

(ii)$\Rightarrow$ (i).  According to Theorem \ref{UWM} we find
that $\{x_k\otimes y_k\}$ is uniformly weakly mixing to zero, for
every bounded sequence $\{y_k\}$ in $Y$. Hence, it is ergodic.
\end{pf}

\section{Applications to $C^*$-dynamical systems}

In this section $\ga$ will be  a $C^*$- algebra with the unity
$\id$. Recall a linear functional $\f\in\ga^*$ is called {\it
positive} if $\f(x^*x)\geq 0$ for every $x\in\ga$. A positive
functional $\f$ is said to be a {\it state} if $\f(\id)=1$. By
$\cs(\ga)$  we denote the set of all states on $\ga$. A linear
operator $T:\ga\mapsto\ga$ is called {\it positive} if $Tx\geq 0$
whenever $x\geq 0$. By $M_n(\ga)$ we denote the set of all
$n\times n$-matrices $a=(a_{ij})$ with entries $a_{ij}$ in $\ga$.
A linear mapping $T:\ga\mapsto\ga$ is called {\it  completely
positive} if the linear operator $T_n:M_n(\ga)\mapsto M_n(\ga)$
given by $T_n(a_{ij})=(T(a_{ij}))$ is positive for all $n\in\bn$.
A completely positive map $T:\ga\mapsto\ga$ with $T\id=\id$ is
called a {\it unital completely positive (ucp)} map.  A pair
$(\ga,T)$ consisting of a $C^*$-algebra $\ga$ and a ucp map
$T:\ga\mapsto\ga$ is called {\it a $C^*$-dynamical system}. Let
$\gb$ be another $C^*$-algebra with unit.  A completion of the
algebraic tensor product $\ga\odot\gb$ with respect to the minimal
$C^*$-tensor norm on $\ga\odot\gb$ is denoted by $\ga\otimes\gb$,
and it would be also a $C^*$-algebra with a unit (see, \cite{T}).
It is known \cite{T} that if $(\ga,T)$ and $(\gb,H)$ are two
$C^*$-dynamical systems, then $(\ga\otimes\gb,T\otimes H)$ is also
$C^*$-dynamical system. Since a mapping $T\otimes
H:\ga\otimes\gb\mapsto\ga\otimes\gb$ given by $(T\otimes
H)(x\otimes y)=Tx\otimes Hy$ is a ucp map.

Let $(\ga,T)$ be a $C^*$–dynamical system, and $\gb$ be a subspace
of $\ga$. Let $E:\ga\to\gb $ be a norm-one projection, i.e. $E^2=E$.
In \cite{F} (see also \cite{AM,FM,MT}) it has been introduced the
following notations

\begin{defin} A $C^*$–dynamical system  $(\ga,T)$ is said to be  \\ \label{smx}
\begin{itemize} \item[(i)] unique $E$--ergodic if
\begin{equation} \label{mp1}
\lim_{n\to\infty}\frac{1}{n}\sum_{k=1}^{n}\f(T^{k}(x))=\f(E(x))\,,\quad
x\in\ga\,,\f\in\cs(\ga)\,.
\end{equation}
\item[(ii)] unique $E$--weakly mixing if
\begin{equation} \label{mpp1}
\lim_{n\to\infty}\frac{1}{n}\sum_{k=0}^{n-1}\big|\f(T^{k}(x))-\f(E(x))\big|=0\,,\quad
x\in\ga\,,\f\in\cs(\ga)\,.
\end{equation}
\end{itemize} \end{defin}

It can readily seen (cf. \cite{AM,FM}) that the map $E$ below is a
norm one projection onto the fixed point subspace
$\ga^{T}=\{x\in\ga: Tx=x\}$. Therefore, in what follows we denote it
by $E_T$. In \cite{AD} (see also \cite{AM}), (i) is called {\it
unique ergodicity w.r.t. the fixed point subalgebra}, whereas (ii)
is called in \cite{FM} {\it $E$--strictly weak mixing}. In addition,
when $E=\w(\,\cdot\,)\id$ (i.e. when there is a unique invariant
state for $T$), (i) is the well--known {\it unique ergodicity}, and
(ii) is called {\it strict (unique) weak mixing} in \cite{MT}. Note
that in \cite{Av} relations between unique ergodicity, minimality
and weak mixing was studied.

By using the Jordan decomposition of bounded linear functionals (cf.
\cite{T}), one can replace $\cs(\ga)$ with $\ga^*$ in Definition
\ref{smx}.

Note that in  \cite{FM,M2} it has been shown that the free shift on
the reduced amalgamated free product $C^{*}$--algebra, and
length--preserving automorphisms of the reduced $C^*$--algebra of
$RD$-group for the length--function, including the free shift on the
free group on infinitely many generators are enjoy unique $E$-mixing
property. Such class of dynamical systems first time was defined and
studied in \cite{AD}. Note that in \cite{FM1} more other complicated
unique $E$-ergodic and unique mixing $C^*$-dynamical systems arising
from free probability have been studied. Note that in \cite{Di}
sufficient and necessary conditions for ergodicity in terms of
joinings are studied.

In this section we are going to apply the results of the previous
section to the given notions.

\begin{thm}\label{mix-a} Let $(\ga,T)$, $(\gb,H)$ be two $C^*$-dynamical systems,
and assume that
 $(\ga\otimes\gb)^*=\ga^*\otimes\gb^*$ is satisfied. Then the following assertions are equivalent:
\begin{itemize}
\item[(i)] The $C^*$-dynamical system $(\ga\otimes\gb,T\otimes H)$
is unique $E_{T\otimes H}$-weak mixing; \item[(ii)] $(\ga,T)$ and
$(\gb,H)$ are unique $E_T$ and $E_H$ weak mixing, respectively.
\end{itemize}

\end{thm}

\begin{pf} (i)$\Rightarrow$(ii) According
to the condition for every an arbitrary functional $\p\in\ga^*$ and
$\v\in \cs(\gb)$, one finds
\begin{eqnarray}\label{mix11}
0&=&
\lim_{n\to\infty}\frac{1}{n}\sum_{k=0}^{n-1}\big|\p\otimes\v(T^k\otimes
H^k(x\otimes \id))-\p\otimes \v (E_{T\otimes
H}(x\otimes\id))\big|\nonumber\\
&=&\lim_{n\to\infty}\frac{1}{n}\sum_{k=0}^{n-1}|\p(T^k(x))-\p\otimes
\v (E_{T\otimes H}(x\otimes\id))|,
\end{eqnarray}
hence
$$
\lim_{n\to\infty}\frac{1}{n}\sum_{k=0}^{n-1}T^k(x)
$$
weak converges, and its limit we denote by $E_T$. Consequently, from
\eqref{mix11} one finds $E_{T\otimes
H}(\cdot\otimes\id)=E_T(\cdot)$. Moreover, $(\ga,T)$ is unique
$E_T$-weak mixing. Similarly, we get unique $E_H$-weak mixing of
$(\gb,H)$.

Let us consider the implication (ii)$\Rightarrow$(i). Let $x\in\ga$
and $y\in\gb$. Define two sequences as follows
\begin{equation}\label{seq}
x_k=T^k(x)-E_T(x), \ \ \ y_k=H^k(y)-E_H(y), \ \ \ k\in\bn.
\end{equation}
Then one can see that the sequences are weakly mixing. Hence,
Theorem \ref{WM} implies that the sequence $\{x_k\otimes y_k\}$ is
weakly mixing as well. This means that for every $\w\in
(\ga\otimes\gb)^*$ one has
\begin{eqnarray}\label{uum1}
\lim_{n\to\infty}\frac{1}{n}\sum_{k=1}^n&&\big|\w(T^k(x)\otimes H^k(y))-\w(T^k(x)\otimes E_H(y))\nonumber\\
&&- \w(E_T(x)\otimes H^k(y))+\w(E_T(x)\otimes E_H(y))\big|=0
\end{eqnarray}

Now define two functionals $\w_1$ and $\w_2$ on $\ga$ and $\gb$,
respectively, as follows:
\begin{equation}\label{ww}
\w_1(\cdot)=\w(\cdot\otimes E_H(y)) \ \ \ \
\w_2(\cdot)=\w(E_T(x)\otimes \cdot),
\end{equation}
here $E_T(x)$ and $E_H(y)$ are fixed. Then according to weak mixing
condition (see (ii)) one has
\begin{eqnarray}\label{eq12}
\lim_{n\to\infty}\frac{1}{n}\sum_{k=1}^{n}\big|\w_1(T^k(x))-\w_1(E_T(x))\big|=0, \\
\label{eq13}
\lim_{n\to\infty}\frac{1}{n}\sum_{k=1}^{n}\big|\w_2(H^k(y))-\w_2(E_H(y))\big|=0.
\end{eqnarray}
The last relations \eqref{eq12},\eqref{eq13} with \eqref{ww} mean
that
\begin{eqnarray}\label{uum2}
\lim_{n\to\infty}\frac{1}{n}\sum_{k=1}^{n}\big|\w(T^k(x)\otimes E_H(y))-\w(E_T(x)\otimes E_H(y))\big|=0, \\
\label{uum3}
\lim_{n\to\infty}\frac{1}{n}\sum_{k=1}^{n}\big|\w(E_T(x)\otimes
H^k(y))-\w(E_T(x)\otimes E_H(y))\big|=0.
\end{eqnarray}

The inequality
\begin{eqnarray*}
|\w(T^k\otimes H^k(x\otimes
y))&-&\w(E_T(x)\otimes E_H(y))| \nonumber \\
&\leq&\bigg|\w(T^k(x)\otimes H^k(y))-\w(T^k(x)\otimes E_H(y))\nonumber\\
&&-
\w(E_T(x)\otimes H^k(y))+\w(E_T(x)\otimes E_H(y))\bigg|\\[2mm]
&&+\big|\w(T^k(x)\otimes E_H(y))-\w(E_T(x)\otimes E_H(y))\big|\\[2mm]
&&+\big|\w(E_T(x)\otimes H^k(y))-\w(E_T(x)\otimes E_H(y))\big|
\end{eqnarray*}
with \eqref{uum1},\eqref{uum2} and \eqref{uum3} imply that
\begin{eqnarray}\label{eq21}
\lim_{n\to\infty}\frac{1}{n}\sum_{k=1}^{n}\big|\w(T^k\otimes
H^k(x\otimes y))-\w(E_T\otimes E_H(x\otimes y))\big|=0.
\end{eqnarray}

The norm-denseness of the elements $\sum_{i=1}^m x_i\otimes y_i$
in $\ga\otimes\gb$ with \eqref{eq21} yields
\begin{eqnarray*}
\lim_{n\to\infty}\frac{1}{n}\sum_{k=1}^{n}\big|\w(T^k\otimes
H^k({\mathbf{z}}))-\w(E_T\otimes E_H({\mathbf{z}}))\big|=0.
\end{eqnarray*}
for arbitrary ${\mathbf{z}}\in\ga\otimes\gb$. So,
$(\ga\otimes\gb,T\otimes H)$ is unique $E_T\otimes E_H$ -weak
mixing.
\end{pf}

\begin{cor}\label{mix-EE} Let $(\ga,T)$ and
$(\gb,H)$ be unique $E_T$ and $E_H$-weak mixing, respectively. Then
one has $E_{T\otimes H}=E_T\otimes E_H$.
\end{cor}

{\bf Remark.} Note that in \cite{L,Wa} certain spectral conditions
of tensor product of dynamical systems defined on von Neumann
algebras were studied. We have to stress that in those papers,
dynamical systems have faithful normal invariant states. For such
weak mixing dynamical systems the condition $E_{T\otimes
H}=E_T\otimes E_H$ is proved as well.

\begin{ex} Now let us provide an example of $C^*$-dynamical system,
which does not have any invariant faithful state, but one has
$E_{T\otimes H}=E_T\otimes E_H$.

Let $\ga=\bc^2$ and $\gb=\bc^3$ and
$$T=
\left(
\begin{array}{ll}
\frac{1}{2} & \frac{1}{2}\\[2mm]
0 & 1
\end{array}
\right), \ \ \ H= \left(
\begin{array}{lll}
1&  0 & 0\\
0& 1 & 0 \\
0 & \frac{1}{2} & \frac{1}{2}
\end{array}
\right).
$$

 It is clear that
\begin{eqnarray*}
&&\ga^T=\{(x,x):\ \ x\in\bc\}, \\
&&\gb^H=\{(x,y,y): \ \ x,y\in\bc\}.
\end{eqnarray*}

One can check that all invariant states for $H$ have the following
form:
$$
(p,q,0), \ \ p,q\geq 0, \ \ p+q=1,
$$
which is not faithful.

 Direct calculations show that
\begin{eqnarray*}
\lim_{n\to\infty}T^n(x,y)=E_T(x,y), \ \
\lim_{n\to\infty}H^n(x,y,z)=E_H(x,y,z),
\end{eqnarray*}
 which mean that $T$ and $H$ are unique $E_T$ and $E_H$ weak mixing,
respectively. Here
$$
E_T(x,y)=(y,y), \ \ \ E_H(x,y,z)=(x,y,y).
$$

Now let us calculate $(\ga\otimes\gb)^{T\otimes H}$. To do it, one
can see that
$$T\otimes H=
\frac{1}{2}\left(
\begin{array}{ll}
H & H\\[2mm]
0 & 2H
\end{array}
\right)
$$

Denote $\xb=(x_1,x_2,x_3),\yb=(y_1,y_2,y_3)$. Then from $T\otimes
H(\xb,\yb)=(\xb,\yb)$ we find
$$
\frac{1}{2}H(\xb+\yb)=\xb, \ \ \ H\yb=\yb.
$$
A simple algebra shows that $\xb=\yb$. Consequently, we have
$$
(\ga\otimes\gb)^{T\otimes H}=\{(x_1,x_2,x_2,x_1,x_2,x_2): \ \
x_1,x_2\in\bc\}
$$
which yields that $(\ga\otimes\gb)^{T\otimes H}=\ga^T\otimes \gb^H$.
This implies that $E_{T\otimes H}=E_T\otimes E_H$.

Moreover, by the same argument we may show that the equality
$E_{H\otimes H}=E_H\otimes E_H$ holds as well.

\end{ex}

{\bf Remark.} The proved theorem extends some results of
\cite{M,M2}. We note that in \cite{Av,L,Wa} similar results were
proved for weak mixing dynamical systems defined over von Neumann
algebras.

Note that some examples of $C^*$-algebras which satisfy the
condition $(\ga\otimes\gb)^*=\ga^*\otimes\gb^*$ can be found in
\cite{M2} (see also \cite{R}).

\begin{thm}\label{mix-c} Let $(\ga,T)$ be a  $C^*$-dynamical systems.
Then for the  following assertions
\begin{itemize}
\item[(i)] $(\ga,T)$ is unique $E_T$-weak mixing;
\item[(ii)] for every $(\gb,H)$ - unique $E_H$-ergodic $C^*$-dynamical
system  with $E_{T\otimes H}=E_T\otimes E_H$ and
$\ga^*\otimes\gb^*=(\ga\otimes\gb)^*$, the $C^*$-dynamical system
$(\ga\otimes\gb,T\otimes H)$ is unique $E_T\otimes E_H$-ergodic;
\end{itemize}
the implication (i)$\Rightarrow$(ii) holds true.
\end{thm}

\begin{pf} Let $(\gb,,H)$ be a $C^*$-dynamical system as in (ii).
Now take arbitrary elements $x\in\ga$ and $y\in\gb$, and consider
the corresponding sequences $\{x_k\}$ and $\{y_k\}$ given by
\eqref{seq}. Then due to the condition $\{x_k\}$ is weak mixing and
$\{y_k\}$ is weak ergodic. Hence, Theorem \ref{UWMM} yields that
$\{x_k\otimes y_k\}$ is weak ergodic, which means for every $\w\in
(\ga\otimes\gb)^*$ one has
\begin{eqnarray}\label{er1}
&&\lim_{n\to\infty}\frac{1}{n}\sum_{k=1}^n\big(\w(T^k(x)\otimes H^k(y))-\w(T^k(x)\otimes E_H(y))\nonumber\\
&&- \w(E_T(x)\otimes H^k(y))+\w(E_T(x)\otimes E_H(y))\big)=0
\end{eqnarray}

Using similar arguments as in the proof of Theorem \ref{mix-a} we
find
\begin{eqnarray}\label{er2}
\lim_{n\to\infty}\frac{1}{n}\sum_{k=1}^{n}\big|\w(T^k(x)\otimes E_H(y))-\w(E_T(x)\otimes E_H(y))\big|=0, \\
\label{er3}
\lim_{n\to\infty}\frac{1}{n}\sum_{k=1}^{n}\big(\w(E_T(x)\otimes
H^k(y))-\w(E_T(x)\otimes E_H(y))\big)=0.
\end{eqnarray}

From
\begin{eqnarray*}
\bigg|\frac{1}{n}\sum_{k=1}^{n}\big(\w(T^k\otimes H^k(x\otimes
y))&-&\w(E_T(x)\otimes E_H(y))\big)\bigg| \nonumber \\
&\leq&\bigg|\frac{1}{n}\sum_{k=1}^{n}\big(\w(T^k(x)\otimes H^k(y))-\w(T^k(x)\otimes E_H(y))\nonumber\\
&&-
\w(E_T(x)\otimes H^k(y))+\w(E_T(x)\otimes E_H(y))\big)\bigg|\\[2mm]
&&+\frac{1}{n}\sum_{k=1}^{n}\big|\w(T^k(x)\otimes E_H(y))-\w(E_T(x)\otimes E_H(y))\big|\\[2mm]
&&+\bigg|\frac{1}{n}\sum_{k=1}^{n}\big(\w(E_T(x)\otimes
H^k(y))-\w(E_T(x)\otimes E_H(y))\big)\bigg|
\end{eqnarray*}
and  \eqref{er1}-\eqref{er3} we obtain
\begin{eqnarray*}
\lim_{n\to\infty}\frac{1}{n}\sum_{k=1}^{n}\big(\w(T^k\otimes
H^k(x\otimes y))&-&\w(E_T\otimes E_H(x\otimes y))\big)=0.
\end{eqnarray*}

Finally, the density argument shows that $(\ga\otimes\gb,T\otimes
H)$ is unique $E_T\otimes E_H$ -ergodic.
\end{pf}

{\bf Remark.} We note that all the results of this section extends
the results of \cite{M,M2} to uniquely $E$-ergodic and uniquely
$E$-weak mixing.

{\bf Remark.} We have to stress that the unique ergodicity $T\otimes
H$ does not imply unique weak mixing of $T$. Indeed, let us consider
the following examples.

\begin{ex} Let $\ga=\bc^2$ and
$$T=
\left(
\begin{array}{ll}
0 & 1\\
1 & 0
\end{array}
\right).
$$
It is clear that $\ga^T=\bc\id$, so $T$ is ergodic, i.e.
$$
\lim_{n\to\infty}\frac{1}{n}\sum_{k=1}^{n}T^k(x,y)=\frac{x+y}{2}(1,1).
\ \ x,y\in\bc
$$
From the equality
$$
\bigg|T^k(x,y)-\frac{x+y}{2}(1,1)\bigg|=\bigg|\frac{x-y}{2}\bigg|
$$
we infer that $T$ is not unique weak mixing.

On the other hand, the equality
$$
(\ga\otimes \ga)^{T\otimes T}=\{(x,y,y,x): \ x,y\in\bc\},
$$
implies  unique $E_{T\otimes T}$-ergodicity of $T\otimes T$.
\end{ex}

\begin{ex} Let $\ga=\bc^3$ and $\gb=\bc^2$. Consider the a
mapping $P:\ga\to\ga$ given by
\begin{equation}\label{3P1}
P(x,y,z)=(y,x,uy+vz),
\end{equation}
where $u,v>0$ and $u+v=1$. It is clear that $P$ is positive and
unital. Direct calculations show that $\ga^P=\bc\id$, which means
$P$ is uniquely ergodic.

Now consider the mapping $P\otimes T$, where $T$ is defined as
above. One can see that such a mapping acts as follows
$$P\otimes T({\mathbf{x}},{\mathbf{y}})=(P{\mathbf{y}},P{\mathbf{x}})
$$
where $\mathbf{x},\mathbf{y}\in\ga$. Hence, we find
$$
(\ga\otimes\gb)^{P\otimes T}=\big\{({\mathbf{x}},P{\mathbf{x}}): \
{\mathbf{x}}\in\ga^{P^2}\big\}.
$$

Therefore, from \eqref{3P1} one immediately gets
\begin{equation}\label{3P2}
P^2(x,y,z)=(x,y,ux+uvy+v^2z).
\end{equation}
Thus, we find
$$
\ga^{P^2}=\bigg\{\bigg(x,y,\frac{x+vy}{1+v}\bigg): \
x,y\in\bc\bigg\}.
$$

On the other hand, we have $\ga^P\otimes \gb^T=\bc\id$, which means
$(\ga\otimes\gb)^{P\otimes T}\neq  \ga^P\otimes \gb^T$.

Similarly reasoning as in Example 3.2 we can show that $P\otimes T$
is uniquely $E_{P\otimes T}$-ergodic.

Note that, from the provided examples we infer the importance of condition
$E_{T\otimes H}=E_T\otimes E_H$.

\end{ex}



\begin{thebibliography}{9999}

\bibitem{ALW} Aaronson, J., Lin, M., Weiss, B. \textit{Mixing
properties of Markov operators and ergodic transformations, and
ergodicity of Cartesian products.} A collection of invited papers
on ergodic theory. Israel J. Math. {\bf 33} (1979), no. 3-4,
198--224.

\bibitem{AD} Abadie B., Dykema K. {\it Unique ergodicity of free shifts and some other automorphisms of
    $C^*$--algebras}, J. Operator Theory, {\bf 61}(2009), 279--294.

\bibitem{AM} Accardi, L., Mukhamedov F., \textit{A note on noncommutative unique
ergodicity and weighted means}, Linear Algebra and Appl.{\bf
430}(2009) 782--790 ({\tt arXiv:0803.0073}).

\bibitem{Av} Avitzour D. {\it Noncommutative topological dynamical systems, II}, Trans. Amer. Math. Soc.
{\bf 282}(1984), 121--135.




\bibitem{BB} Berend D., Bergelson V., \textit{Mixing sequences in Hilbert spaces}, Proc. Amer. Math. Soc.
{\bf 98}(1986), 239--246.

\bibitem{BH} Blum J.R., Hanson D.L., \textit{On the mean ergodic theorem for
subsequences}, Bull.Amer. Math. Soc. \textbf{66} (1960), 308--311.


\bibitem{D} Diestel, J. \textit{Sequences and series in Banach spaces}, Graduate texts in Math. 92,
Springer-Verlag, 1984.

\bibitem{Di} Duvenhage, R. \textit{Joinings of W*-dynamical systems}, Jour. Math. Anal. Appl.
{\bf 343}(2008), 175--181.


\bibitem{F} Fidaleo, F., \textit{On strong ergodoc properties of quantum dynamical systems}, Infinite Dimen.
Anal. Quantum Probab. Related Topics {\bf 12}(2009), 551--556.


\bibitem{FM} Fidaleo, F., Mukhamedov F., \textit{Strict weak mixing of some $C^*$-dynamical systems
based on free shifts}, Jour. Math. Anal. Appl. {\bf 336}(2007),
180--187.

\bibitem{FM1} Fidaleo, F., Mukhamedov F., \textit{Ergodic properties of Bogoliubov
automorphisms in free probability}, Inf. Dim. Anal. Quantum Probab.
and Related Topics, {\bf 13}(2010) 393--411

\bibitem{JL} Jones L.K., Lin, M., \textit{Ergodic theorems of weak mixing type.}
Proc. Amer. Math. Soc. {\bf 57} (1976), 50--52.

\bibitem{K} Krengel H.O., \textit{Ergodic Theorems}, Walter de Gruyter,
Berlin-New York, 1985.

\bibitem{L} {\L}uczak A.
{\it Eigenvalues and eigenspaces of quantum dynamical systems and
their tensor products}, J. Math. Anal. Appl. {\bf 221} (1998),
13--32.

\bibitem{M} Mukhamedov F., \textit{On strictly weakly mixing $C^*$-dynamical
systems}, Funct. Anal. Appl. {\bf 27}(2007), 311--313.

\bibitem{M2} Mukhamedov F. \textit{On strictly weak mixing C*-dynamical systems and  a
weighted ergodic theorem}, Studia Sci. Math. Hungarica {\bf
47}(2010), 155--174.

\bibitem{MT} Mukhamedov F., Temir S. \textit{A few remarks on mixing properties of
$C^*$-dynamical systems}, Rocky Mount. J. Math. {\bf 37}(2007),
1685--1703.

\bibitem{NSZ} Nicolescu, C., Str\"oh, A., Zsid\'o, L., \textit{Noncommutative
extensions of classical and multiple recurrence theorems,} J.
Operator Theory, {\bf 50}(2003), 3-52.


\bibitem{R} Ruan R.A. \textit{Introduction to Tensor Products of Banach Spaces},
Springer, London--Berlin--Heidelberg, 2002.

\bibitem{T} Takesaki, M., \textit{Theory of Operator algebras, I}, Springer,
Berlin--Heidelberg--New York, 1979.

\bibitem{Wa} Watanabe, S. {\it Asymptotic behavior and eigenvalues of
dybamical semi-groups on operator algebars}, J. Math. Anal. Appl.
{\bf 86} (1982), 411--424.

\bibitem{Z}  Zsid\'o L. Weak mixing properties of
vector sequences, In book: Dritschel, M.A. (ed.), {\it The extended
field of operator theory}. Containing lectures of the 15th
international workshop on operator theory and its applications,
IWOTA 2004, Newcastle, UK, July 12--16, 2004. Basel: Birkhauser.
Operator Theory: Advances and Applications 171, 361--388 (2006).



\end{thebibliography}
\end{document}